\documentclass[a4paper,12pt]{article}
\usepackage{amscd}
\usepackage{amsmath,amsfonts,amssymb,amscd}
\usepackage{indentfirst,graphicx,epsfig}
\usepackage{graphicx,psfrag}
\usepackage{ifpdf}
\usepackage{float}

\input{epsf}

\newtheorem{thm}{Theorem}
\newtheorem{pro}{Proposition}

\newtheorem{lem}{Lemma}

\newtheorem{cla}{Claim}

\def\qed{\hfill \nopagebreak\rule{5pt}{8pt}}
\def\pf{\noindent {\it Proof.} }

\begin{document}
\rule{0cm}{1cm}
\begin{center}
{\Large\bf Further hardness results on the\\[3mm] rainbow
vertex-connection number of graphs}
\end{center}

\begin{center}
Lily Chen, Xueliang Li, Huishu Lian\\
Center for Combinatorics, LPMC\\
Nankai University, Tianjin 300071, P. R. China\\
Email:  lily60612@126.com, lxl@nankai.edu.cn, lhs6803@126.com
\end{center}

\begin{abstract}
A vertex-colored graph $G$ is {\it rainbow vertex-connected} if any
pair of vertices in $G$ are connected by a path whose internal
vertices have distinct colors, which was introduced by Krivelevich
and Yuster. The {\it rainbow vertex-connection number} of a
connected graph $G$, denoted by $rvc(G)$, is the smallest number of
colors that are needed in order to make $G$ rainbow
vertex-connected. In a previous paper we showed that it is
NP-Complete to decide whether a given graph $G$ has $rvc(G)=2$. In
this paper we show that for every integer $k\geq 2$, deciding
whether $rvc(G)\leq k$ is NP-Hard. We also show that for any fixed
integer $k\geq 2$, this problem belongs to NP-class, and so it becomes
NP-Complete. \\[3mm]
\end{abstract}

\section{Introduction}

All graphs considered in this paper are simple, finite and
undirected. Undefined terminology and notation can be found in
\cite{BM}.

Let $G$ be a nontrivial connected graph with an edge-coloring $c:
E(G)\rightarrow \{1,2,\cdots,k\}$, $k \in \mathbb{N}$, where
adjacent edges may be colored the same. A path $P$ of $G$ is a
\emph{rainbow path} if no two edges of $P$ are colored the same. The
graph $G$ is called \emph{rainbow-connected} if for any pair of
vertices $u$ and $v$ of $G$, there is a rainbow $u-v$ path. The
minimum number of colors for which there is an edge-coloring of $G$
such that $G$ is rainbow connected is called the \emph{rainbow
connection number}, denoted by $rc(G)$. Clearly, if a graph is
rainbow connected, then it is also connected. Conversely, any
connected graph has a trivial edge-coloring that makes it rainbow
connected, just assign each edge a distinct color. An easy
observation is that if $G$ has $n$ vertices then $rc(G)\leq n-1$,
since one may color the edges of a spanning tree with distinct
colors, and color the remaining edges with one of the colors already
used. It is easy to see that if $H$ is a connected spanning subgraph
of $G$, then $rc(G)\leq rc(H)$. We note the trivial fact that
$rc(G)=1$ if and only if $G$ is a clique, the fact that $rc(G)=n-1$
if and only if $G$ is a tree, and the easy observation that a cycle
with $k\geq 4$ vertices has a rainbow connection number $\lceil
k/2\rceil$. Also notice that $rc(G)\geq diam(G)$, where $diam(G)$ is
the diameter of $G$.

Similar to the concept of rainbow connection number, Krivelevich and
Yuster \cite{KY} proposed the concept of rainbow vertex-connection.
Let $G$ be a nontrivial connected graph with a vertex-coloring $c:
V(G)\rightarrow \{1,2,\cdots,k\}, k \in \mathbb{N}$. A path $P$ of
$G$ is \emph{rainbow vertex-connected} if its internal vertices have
distinct colors. The graph $G$ is {\it rainbow vertex-connected} if
any pair of vertices are connected by a rainbow vertex-connected
path. In particular, if $k$ colors are used, then $G$ is rainbow
$k$-vertex-connected. The {\it rainbow vertex-connection number} of a connected
graph $G$, denoted by $rvc(G)$, is the smallest number of colors
that are needed in order to make $G$ rainbow vertex-connected. An
easy observation is that if $G$ is of order $n$ then $rvc(G)\leq
n-2$, $rvc(G)=0$ if and only if $G$ is a complete graph, and
$rvc(G)=1$ if and only if $diam(G)=2$. Notice that $rvc(G)\geq
diam(G)-1$ with equality if the diameter is $1$ or $2$. For the
rainbow connection number and the rainbow vertex-connection number,
some examples were given to show that there is no upper bound for
one of parameters in terms of the other in \cite{KY}. Krivelevich
and Yuster \cite{KY} proved that if $G$ is a graph with $n$ vertices
and minimum degree $\delta$, then $rvc(G)<11n/\delta$. Li and Shi
used a similar proof technique and greatly improved this bound, see
\cite{LS}.

The computational complexity of rainbow connection number has been
studied extensively. In \cite{CLRTY}, Caro et al. conjectured that
computing $rc(G)$ is an NP-Hard problem, and that even deciding
whether a graph has $rc(G)=2$ is NP-Complete. Later, Chakraborty et
al. confirmed this conjecture in \cite{CFMY}. They also conjectured
that for every integer $k\geq 2$, to decide whether $rc(G)\leq k$ is
NP-Hard. Recently, Ananth and Nasre confirmed the conjecture in
\cite{AN}. Li and Li \cite{SX} showed that for any fixed integer
$k\geq 2$, to decide whether $rc(G)\leq k$ is actually NP-Complete.
For the rainbow vertex-connection number we got a similar complexity
result in \cite{CLS}.
\begin{thm} \cite{CLS} \label{theorem1}
Given a graph $G$, deciding whether $rvc(G)=2$ is NP-Complete. Thus,
computing $rvc(G)$ is NP-Hard.
\end{thm}

As a generalization of the above result, in this paper we will show
the following result:
\begin{thm}\label{theorem2}
For every integer $k\geq 2$, to decide whether $rvc(G)\leq k$ is
NP-Hard. Moreover, for any fixed integer $k\geq 2$, the problem
belongs to NP-class, and therefore it is NP-Complete.
\end{thm}

In order to prove this theorem, we first show that an intermediate
problem called the $k$-subset rainbow vertex-connection problem is
NP-Hard by giving a reduction from the vertex-coloring problem. We
then establish the polynomial-time equivalence of the $k$-subset
rainbow vertex-connection problem and the problem of deciding
whether $rvc(G)\leq k$ for a graph $G$.

\section{Proof of Theorem 2}

We first describe the problem of $k$-subset rainbow
vertex-connection: given a graph $G$ and a set of pairs $P\subseteq
V(G)\times V(G)$, decide whether there is a vertex-coloring of $G$
with $k$ colors such that every pair of vertices $(u,v)\in P$ is
rainbow vertex-connected. Recall that the $k$-vertex-coloring
problem is as follows: given a graph $G$ and an integer $k$, whether
there exists an assignment of at most $k$ colors to the vertices of
$G$ such that no pair of adjacent vertices are colored the same. It
is known that this $k$-vertex-coloring problem is NP-Hard for $k\geq
3$. Now we reduce the $k$-vertex-coloring problem to the $k$-subset
rainbow vertex-connection problem, which shows that the problem of
$k$-subset rainbow vertex-connection is NP-Hard.
\begin{lem}
The problem of $k$-vertex-coloring is polynomially reducible to the
problem of $k$-subset rainbow vertex-connection.
\end{lem}
\pf Let $G=(V,E)$ be an instance of the $k$-vertex-coloring problem,
we construct a graph $\langle G'=(V',E'),P\rangle$ as follows:

For every vertex $v\in V$ we introduce a new vertex $x_v$. We set
$$V'=V\cup \{x_v: \ v\in V \} \ and \ E'=E\cup\{(v,x_v): \ v\in V\}.$$

Now we define the set $P$ as follows: $$P=\{(x_u,x_v):(u,v)\in
E\}.$$

It remains to verify that $G$ is vertex-colorable using $k(\geq 3)$
colors if and only if there is a vertex-coloring of $G'$ with $k$
colors such that every pair of vertices $(x_u,x_v)\in P$ is rainbow
vertex-connected.

Let $c$ be the proper $k$-vertex-coloring of $G$. We define the
vertex-coloring $c'$ of $G'$ by $c'(x_v)=c'(v)=c(v)$. If
$(x_u,x_v)\in P,$ then $(u,v)\in E$, $c(u)\neq c(v),$  and so
$c'(u)\neq c'(v)$, $x_uuvx_v$ is a rainbow vertex-connected path
between $x_u$ and $x_v$.

In the other direction, assume that $c'$ is a $k$-vertex-coloring of
$G'$ such that every pair of vertices $(x_u,x_v)\in P$ is rainbow
vertex-connected. We define the vertex-coloring $c$ of $G$ by
$c(v)=c'(v)$. For every $(u,v)\in E$, $(x_u,x_v)\in P$, since the
rainbow vertex-connected path between $x_u$ and $x_v$ must go
through $u$ and $v$, $c'(u)\neq c'(v)$, and so $c(u)\neq c(v)$, thus
$c$ is the proper $k$-vertex-coloring of $G$. \qed

In the following, we prove that the problem of deciding whether a
graph is $k$-subset rainbow vertex-connection is polynomial-time
equivalent to the problem of deciding whether $rvc(G)\leq k$ for a
graph $G$.
\begin{lem}
The following problems are polynomial-time equivalent: \\
1. Given a graph $G$, decide whether $rvc(G)\leq k$. \\
2. Given a graph $G$ and a set $P\subseteq V(G)\times V(G)$ of pairs
of vertices, decide whether there is a vertex-coloring of $G$ with
$k$ colors such that every pair of vertices $(u,v)\in P$ is rainbow
vertex-connected.
\end{lem}
\pf It is sufficient to demonstrate a reduction from Problem $2$ to
Problem $1$. Let $\langle G=(V,E), P\rangle$ be any instance of
Problem $2$. We construct a graph $G_k=(V_k, E_k)$ such that $G$ is
a subgraph of $G_k$ and $rvc(G_k)\leq k$ if and only if $G$ is
$k$-subset rainbow vertex-connected. We prove the correctness of the
reduction by induction on $k$. For $k=2$ and $k=3$, we give explicit
constructions and show that the reduction is valid. Then we show our
inductive step to get $G_k$ and prove the correctness of the
reduction.

\noindent {\bf Construction of $G_2$:} Let $G_2=(V_2,E_2)$ where the
vertex set $V_2$ is defined as follows:
\begin{eqnarray*}
V_2&=&\{u\} \cup V_2^{(0)}\cup V_2^{(2)}\\
V_2^{(0)}&=&\{v_{i,0}^{(1)}, v_{i,0}^{(2)}:i \in \{1, 2, \cdots,
n\}\}
\cup \{w_{i,j}^{(1)},w_{i,j}^{(2)} :(v_i,v_j)\in(V \times V)\backslash{P}\}\\
V_2^{(2)}&=&\{v_{i,2}: i\in \{1, 2, \cdots n\}\}
\end{eqnarray*}
and the edge set $E_2$ is defined as:
\begin{eqnarray*}
E_2&=& E_2^{(1)}\cup E_2^{(2)}\cup E_2^{(3)}\cup E_2^{(4)}\cup E_2^{(5)}\cup E_2^{(6)}\\
E_2^{(1)}&=&\{(u,x):x \in V_2^{(0)}\}\\
E_2^{(2)}&=&\{(v_{i,0}^{(1)}, v_{i,0}^{(2)}):i\in \{1, 2, \cdots, n\}\}\\
E_2^{(3)}&=&\{(w_{i,j}^{(1)},w_{i,j}^{(2)}) :(v_i,v_j)\in(V \times V)\backslash{P}\}\\
E_2^{(4)}&=&\{(v_{i,2},v_{i,0}^{(1)}),(v_{i,2},v_{i,0}^{(2)}):i\in \{1, 2,\cdots,n\}\}\\
E_2^{(5)}&=&\{(v_{i,2},w_{i,j}^{(1)}),(v_{j,2},w_{i,j}^{(2)}):(v_i,v_j)\in(V\times
V)\backslash{P}\}\\
E_2^{(6)}&=&\{(v_{i,2},v_{j,2}):(v_i,v_j)\in E(G)\}
\end{eqnarray*}

Denote $H_2=G_2[\{v_{i,2}:i \in \{1,2,\ldots,n\}\}]$. Let
$P_2=\{(v_{i,2},v_{j,2}):(v_i,v_j)\in P\}$. The graph $G_2$
satisfies the property that for all $(v_{i,2},v_{j,2})\in P_2$ there
is no path of length $\leq 3$ between $v_{i,2}$ and $v_{j,2}$ in
$G_2\setminus E(H_2)$ and also for all $(v_{i,2},v_{j,2})\notin P_2$
the length of the shortest path between $v_{i,2}$ and $v_{j,2}$ in
$G_2\setminus E(H_2)$ is $3$.

Let $c: V\rightarrow \{1,2\}$ be a $2$-vertex-coloring of $G$ such
that every pair of vertices in $P$ is rainbow vertex-connected .
Define the vertex-coloring $c_2$ of $G_2$ as follows:
\begin{itemize}
\item $c_2(u)=1$.
\item $c_2(v_{i,0}^{(1)})=1$ and $c_2(v_{i,0}^{(2)})=2$ for $i\in \{1, 2, \cdots,
n\}$.

$c_2(w_{i,j}^{(1)})=1$ and $c_2(w_{i,j}^{(2)})=2$, for all
$w_{i,j}^{(\alpha)} \in V_2^{(0)},\alpha \in \{1,2\}$.
\item $c_2(v_{i,2})=c(v_i)$, for $i\in \{1, 2, \cdots, n\}$.
\end{itemize}
It can be easily verified that $rvc(G_2)\leq 2$ if and only if $G$
is $2$-subset rainbow vertex-connected.

\noindent {\bf Construction of $G_3$:} Let $G_3=(V_3,E_3)$ where the
vertex set $V_3$ is defined as follows:
\begin{eqnarray*}
V_3&=&V_3^{(0)} \cup V_3^{(1)}\cup V_3^{(3)}\\
V_3^{(0)}&=&\{v_{i,0}^{(1)}, v_{i,0}^{(2)}:i\in \{1, 2, \cdots, n\}\}\cup \{u_{i,j}^{(1)},u_{i,j}^{(2)}: (v_i,v_j)\in(V \times V)\backslash{P}\}\\
V_3^{(1)}&=&\{v_{i,1}^{(1)}, v_{i,1}^{(2)}:i\in \{1, 2, \cdots, n\}\}\cup \{w_{i,j}^{(1)},w_{i,j}^{(2)} :(v_i,v_j)\in(V \times V)\backslash{P}\}\\
V_3^{(3)}&=&\{v_{i,3}: i\in \{1, 2, \cdots, n\}\}
\end{eqnarray*}
and the edge set $E_3$ is defined as:
\begin{eqnarray*}
E_3&=&E_3^{(1)}\cup E_3^{(2)}\cup E_3^{(3)}\cup E_3^{(4)}\cup E_3^{(5)}\cup E_3^{(6)}\cup E_3^{(7)}\\
E_3^{(1)}&=&\{(x,y):x,y \in V_3^{(0)}\}\\
E_3^{(2)}&=&\{(v_{i,0}^{(\alpha)},v_{i,1}^{(\beta)}):i\in\{1,2,\cdots,n\}, \ \alpha,\beta\in \{1,2\}\} \\
E_3^{(3)}&=&\{(u_{i,j}^{(\alpha)},w_{i,j}^{(\beta)}) :(v_i,v_j)\in(V \times V)\backslash{P},\ \alpha,\beta \in \{1,2\}\} \\
E_3^{(4)}&=&\{(v_{i,1}^{(1)}, v_{i,1}^{(2)}):i\in \{1, 2, \cdots, n\}\}\\
E_3^{(5)}&=&\{(v_{i,3},v_{i,1}^{(1)}),(v_{i,3},v_{i,1}^{(2)}):i\in
\{1,2,\cdots, n\}\}\\
E_3^{(6)}&=&\{(v_{i,3},w_{i,j}^{(1)}),(v_{j,3},w_{i,j}^{(2)}):(v_i,v_j)\in(V\times
V)\backslash{P}\}\\
E_3^{(7)}&=&\{(v_{i,3},v_{j,3}):(v_i,v_j)\in E(G)\}
\end{eqnarray*}

Denote $H_3=G_3[\{v_{i,3}:i \in \{1,2,\ldots,n\}\}]$. Let
$P_3=\{(v_{i,3},v_{j,3}):(v_i,v_j)\in P\}$. The graph $G_3$
satisfies the property that for all $(v_{i,3},v_{j,3})\in P_3$ there
is no path of length $\leq 4$ between $v_{i,3}$ and $v_{j,3}$ in
$G_3\setminus E(H_3)$ and also for all $(v_{i,3},v_{j,3})\notin
 P_3$ the length of the shortest path between
$v_{i,3}$ and $v_{j,3}$ in $G_3\setminus E(H_3)$ is $4$.

Let $c: V\rightarrow \{1,2,3\}$ be a $3$-vertex-coloring of $G$ such
that every pair of vertices in $P$ is rainbow vertex-connected.
Define the vertex-coloring $c_3$ of $G_3$ as follows:
\begin{itemize}
\item $c_3(v_{i,0}^{(1)})=1$ and $c_3(v_{i,0}^{(2)})=2$, for $i\in \{1, 2, \cdots,
n\}$,

$c_3(u_{i,j}^{(1)})=1$ and $c_3(u_{i,j}^{(2)})=2$, for
$u_{i,j}^{(1)},u_{i,j}^{(2)} \in V_3^{(0)}$.
\item $c_3(v_{i,1}^{(1)})=2$ and $c_3(v_{i,1}^{(2)})=3$, for $i\in \{1, 2, \cdots,
n\}$,

$c_3(w_{i,j}^{(1)})=2$ and $c_3(w_{i,j}^{(2)})=3$, for
$w_{i,j}^{(1)}, w_{i,j}^{(2)} \in V_3^{(1)}$.
\item $c_3(v_{i,3})=c(v_i)$, for $i\in \{1, 2, \cdots, n\}$.
\end{itemize}
It can be easily verified that $rvc(G_3)\leq 3$ if and only if $G$
is $3$-subset rainbow vertex-connected.

\noindent {\bf Inductive construction of $G_k$:} Assuming that we
have constructed $G_{k-2}=(V_{k-2},E_{k-2})$, the graph
$G_k=(V_k,E_k)$ is then constructed as follows: Each base vertex
$v_{i,k-2}$ in $V_{k-2}$ is split into the vertices
$v_{i,k-2}^{(1)},v_{i,k-2}^{(2)}$ and edges are added between them.
Any edge of the form $(x,v_{i,k-2})$ is replaced by
$(x,v_{i,k-2}^{(1)})$, $(x,v_{i,k-2}^{(2)})$. After doing this, we
add the vertices $v_{i,k}$ and edges $(v_{i,k},v_{i,k-2}^{(1)})$,
$(v_{i,k},v_{i,k-2}^{(2)})$ for $i\in \{1,2, \cdots, n\}$. Formally
the graph $G_k$ is defined as follows:

When $k$ is even: $V_k=\{u\} \cup V_k^{(0)} \cup V_k^{(2)} \cup
\cdots \cup V_k^{(k)}$, where
\begin{eqnarray*}
\begin{split}
V_k^{(i)}=&V_{k-2}^{(i)}, \ \ \ \   for\ \ \   i=0,2,\cdots,k-4;\\
V_k^{(k-2)}=&\{v_{i,k-2}^{(1)},v_{i,k-2}^{(2)}: i \in \{1,2,\cdots,n\}\};\\
V_k^{(k)}=&\{v_{i,k}:i \in \{1,2,\cdots,n\}\}. \\
\end{split}
\end{eqnarray*}

When $k$ is odd: $V_k=V_k^{(0)} \cup V_k^{(1)} \cup V_k^{(3)}\cup
\cdots \cup V_k^{(k)}$, where
\begin{eqnarray*}
\begin{split}
V_k^{(i)}=&V_{k-2}^{(i)}, \ \ \ \   for\ \ \   i=0,1,3,\cdots,k-4;\\
V_k^{(k-2)}=&\{v_{i,k-2}^{(1)},v_{i,k-2}^{(2)}: i \in \{1,2,\cdots,n\}\};\\
V_k^{(k)}=&\{v_{i,k}:i \in \{1,2,\cdots,n\}\}. \\
\end{split}
\end{eqnarray*}

For all $k\geq 4$, $E_k$ is defined as follows:
\begin{eqnarray*}
\begin{split}
E_k=&E_{k-2}\setminus (E_{k-2}(V_{k-2}^{(k-4)},V_{k-2}^{(k-2)}) \cup E(H_{k-2}))\\
&\cup \{(v_{i,k-2}^{(\alpha)},x):(v_{i,k-2},x)\in E_{k-2}(V_{k-2}^{(k-4)},V_{k-2}^{(k-2)}),i\in \{1, 2, \cdots n\},\alpha \in \{1,2\}\}\\
&\cup \{(v_{i,k-2}^{(1)},v_{i,k-2}^{(2)}): i\in
\{1,2,\cdots,n\}\}\\
&\cup \{(v_{i,k},v_{i,k-2}^{(\alpha)}):i\in \{1,2, \cdots,
n\},\alpha \in \{1,2\}\}\cup E(H_k)\\
\end{split}
\end{eqnarray*}
where $E(H_l)=\{(v_{i,l},v_{j,l}): (v_i,v_j)\in E(G)\}$  and
$E_{k-2}(V_{k-2}^{(k-4)},V_{k-2}^{(k-2)})=\{(u,v):u\in
V_{k-2}^{(k-4)}, v\in V_{k-2}^{(k-2)}\}$.

Let $P_k=\{(v_{i,k},v_{j,k}):(v_i,v_j)\in P\}$. Then we show that
the graph $G_k$ satisfies the following properties as claims:
\begin{cla}\label{Claim1}
For any $(v_{i,k},v_{j,k})\in P_k$, there is no path of length less
than $k+2$ between $v_{i,k}$ and $v_{j,k}$ in $G_k \setminus
E(H_k)$.
\end{cla}
\pf It has been shown that the assertion is true for $G_2$ and
$G_3$. Assume that the assertion is true for $G_{k-2}$. Let
$(v_i,v_j)\in P$, then $(v_{i,k-2},v_{j,k-2})\in P_{k-2}$, and hence
by induction, there is no path of length less than $k$ between
$v_{i,k-2}$ and $v_{j,k-2}$ in $G_{k-2} \setminus E(H_{k-2})$. By
the construction of $G_k$, we do not shorten the paths between any
two vertices, so the paths from $v_{i,k-2}^{(\alpha)}$ to
$v_{j,k-2}^{(\beta)}$ will still be of length at least $k$ for
$\alpha, \beta \in \{1,2\}$. Consider the graph $G_k \setminus
E(H_k)$. Since the neighbors of the vertex $v_{i,k}$ are only
$v_{i,k}^{(1)}, v_{i,k}^{(2)}$, the path between $v_{i,k}$ and
$v_{j,k}$ must be $v_{i,k}v_{i,k-2}^{(\alpha)}\ldots
v_{j,k-2}^{(\beta)}v_{j,k}$ for $\alpha=1\ or\ 2, \beta=1\ or\ 2$,
thus their lengths are at least $k+2$.\qed
\begin{cla}
For any $(v_{i,k},v_{j,k})\notin P_k$, the shortest path between
$v_{i,k}$ and $v_{j,k}$ is of length $k+1$ in $G_k \setminus
E(H_k)$.
\end{cla}
\pf It has been shown that the assertion is true for $G_2$ and
$G_3$. Suppose that the assertion is true for $G_{k-2}$. Let
$(v_i,v_j)\notin P$, then $(v_{i,k-2},v_{j,k-2})\notin P$, and hence
by induction, the shortest path between $v_{i,k-2}$ and $v_{j,k-2}$
is of length $k-1$ in $G_{k-2} \setminus E(H_{k-2})$. By the
construction of $G_k$, we do not shorten the paths between any two
vertices, so the shortest path between $v_{i,k-2}^{(\alpha)}$ and
$v_{j,k-2}^{(\beta)}$ will still be of length  $k-1$ for $\alpha,
\beta \in \{1,2\}$. Consider the graph $G_k \setminus E(H_k)$. Since
the neighbors of the vertex $v_{i,k}$ are only $v_{i,k}^{(1)},
v_{i,k}^{(2)}$, the shortest path between $v_{i,k}$ and $v_{i,k}$
must be $v_{i,k}v_{i,k-2}^{(\alpha)}\cdots
v_{j,k-2}^{(\beta)}v_{j,k-2}$ for $\alpha=1\ or\ 2, \beta=1\ or\ 2$,
thus the length of the path is $k+1$.\qed
\begin{cla}
$G$ is $k$-subset rainbow vertex-connected if and only if $G_k$ is
$k$-rainbow vertex-connected.
\end{cla}
\pf  Denote $H_k=G_k[\{v_{i,k}: i \in \{1,2,\cdots,n\}\}]$. It can
be seen that $H_k$ is isomorphic to $G$.

If $G_k$ is $k$-rainbow vertex-connected, let $c_k:
V(G_k)\rightarrow \{1,2,\cdots, k\}$ be a vertex-coloring of $G_k$
with $k$ colors such that every pair of vertices in $G_k$ is rainbow
vertex-connected. We define the vertex-coloring $c$ of $G$ as
follows: $c(v_i)=c_k(v_{i,k})$ for $i \in \{1,2,\cdots,n\}$. If
$(v_i,v_j)\in P$, then $(v_{i,k},v_{j,k})\in P_k$. By Claim
\ref{Claim1}, there is no path between $v_{i,k}$ and $v_{j,k}$ with
length less than $k+2$ in $G_k\setminus E(H_k)$. Hence the entire
rainbow vertex-connected path between $v_{i,k}$ and $v_{j,k}$ must
lie in $H_k$ itself. Correspondingly, there is a rainbow
vertex-connected path between $v_i$ and $v_j$ in $G$. Thus, $G$ is
$k$-subset rainbow vertex-connected.

In the other direction, if $G$ is $k$-subset rainbow
vertex-connected, let $c: V(G)\rightarrow \{1,2,\cdots, k\}$ be a
vertex-coloring of $G$ with $k$ colors such that every pair of
vertices in $P$ is rainbow vertex-connected. We define the
vertex-coloring $c_k$ of $G_k$ by induction. We have given the
vertex-colorings $c_2, c_3$ of $G_2, G_3$. Assume that $c_{k-2}:
V(G_{k-2})\rightarrow \{1,2,\cdots, k-2\}$ is a vertex-coloring of
$G_{k-2}$ such that $G_{k-2}$ is rainbow vertex-connected. We define
the vertex-coloring $c_k$ of $G_k$ as follows:

When $k$ is even:
\begin{itemize}
\item $c_k(u)=k-1$.
\item $c_k(v)=c_{k-2}(v)$, for $v\in V_k^{(0)}\cup V_k^{(2)}\cup \cdots \cup V_k^{(k-4)}$.
\item $c_k(v_{i,k-2}^{(1)})=k-1, c_k(v_{i,k-2}^{(2)})=k$, for $i\in \{1, 2, \cdots, n\}$.
\item $c_k(v_{i,k})=c(v_i)$, for $i\in \{1, 2, \cdots, n\}$.
\end{itemize}

When $k$ is odd:
\begin{itemize}
\item
$c_k(v_{i,0}^{(1)})=c_{k-2}(v_{i,0}^{(1)}), c_k(v_{i,0}^{(2)})=k-1$,
for $i \in \{1,2,\cdots,n\}$.

$c_k(u_{i,j}^{(1)})=c_{k-2}(u_{i,j}^{(1)}), c_k(u_{i,j}^{(2)})=k-1$
for $u_{i,j}^{(1)},u_{i,j}^{(2)} \in V_k^{(0)}$.
\item $c_k(v)=c_{k-2}(v)$, for $v\in V_k^{(1)}\cup V_k^{(3)}\cup \cdots \cup V_k^{(k-4)}$.
\item $c_k(v_{i,k-2}^{(1)})=k-1, c_k(v_{i,k-2}^{(2)})=k$, for $i\in \{1, 2, \cdots, n\}$.
\item $c_k(v_{i,k})=c(v_i)$, for $i\in \{1, 2, \cdots, n\}$.
\end{itemize}

\begin{pro}
The vertex-coloring $c_k$ of $G_k$ defined above makes $G_k$ rainbow
vertex-connected.
\end{pro}

\pf Let $v, w\in V_k,$ we now show that $v, w$ are rainbow
vertex-connected in $G_k$.

{\bf Case 1.} $k$ is even.

By the vertex-coloring $c_k$, we have $c_k(v_{i,j}^{(1)})=j+1,$
$c_k(v_{i,j}^{(2)})=j+2,$  $c_k(u)=k-1$ and $c_k(v_{i,k})=c(v_i)$
for $i\in\{1,2,\cdots,n\},j\in \{0,2,\cdots,k-2\}$.

{\bf Subcase 1.1.} $v\in V_k^{(p)}$, $w\in V_k^{(q)}$, where $p,q\in
\{0,2,\cdots, k-2\}$.

If $v=v_{i,p}^{(\alpha)}$, $w=v_{j,q}^{(\beta)}$ for
$\alpha,\beta\in \{1,2\},$ then
$vv_{i,p-2}^{(1)}v_{i,p-4}^{(1)}\cdots
v_{i,0}^{(1)}uv_{j,0}^{(2)}\cdots v_{j,q-2}^{(2)}w$ is the rainbow
vertex-connected path between $v$ and $w$.

If $v=v_{i_1,p}^{(\alpha)}$, $w=w_{i,j}^{(\beta)}$ for
$\alpha,\beta\in \{1,2\},$ then
$vv_{i_1,p-2}^{(1)}v_{i_1,p-4}^{(1)}\cdots v_{i_1,0}^{(1)}uw$ is the
rainbow vertex-connected path between $v$ and $w$.

If $v=w_{i_1,j_1}^{(\alpha)}$, $w=w_{i_2,j_2}^{(\beta)}$ for
$\alpha,\beta\in \{1,2\},$ then $vuw$ is the rainbow
vertex-connected path between $v$ and $w$.

{\bf Subcase 1.2.} $v=v_{i,k}$, $w\in V_k^{(q)}$, where $q\in
\{0,2,\cdots, k-2\}$.

If $w=v_{j,q}^{(\alpha)}$ for $\alpha\in \{1,2\},$ then
$vv_{i,k-2}^{(2)}v_{i,k-4}^{(1)}\cdots
v_{i,0}^{(1)}uv_{j,0}^{(2)}\cdots v_{j,q-2}^{(2)}w$ is the rainbow
vertex-connected path between $v$ and $w$.

If $w=w_{i,j}^{(\alpha)}$ for $\alpha \in \{1,2\},$ then
$vv_{i,k-2}^{(2)}v_{i,k-4}^{(1)}\cdots v_{i,0}^{(1)}uw$ is the
rainbow vertex-connected path between $v$ and $w$.

{\bf Subcase 1.3.} $v=v_{i,k}$, $w=v_{j,k}$.

If $(v_{i,k},v_{j,k})\in P_k$, then $(v_i,v_j)\in P$. By the
vertex-coloring $c$ of $G$, there is a rainbow vertex-connected path
between $v_i$ and $v_j$ in $G$. Correspondingly, since
$c_k(v_{i,k})=c(v_i)$, there is a rainbow vertex-connected path
between $v_{i,k}$ and $v_{j,k}$ in $G_k$.

If $(v_{i,k},v_{j,k})\notin P_k$, then
$v_{i,k}v_{i,k-2}^{(1)}v_{i,k-4}^{(1)}\cdots
v_{i,2}^{(1)}w_{i,j}^{(1)}w_{i,j}^{(2)}v_{j,2}^{(2)}\cdots
v_{j,k-2}^{(2)}v_{j,k}$ is the rainbow vertex-connected path between
$v_{i,k}$ and $v_{j,k}$.

{\bf Case 2.} $k$ is odd.

By the vertex-coloring $c_k$, we have

$c_k(v_{i,j}^{(1)})=j+1,$ $c_k(v_{i,j}^{(2)})=j+2,$ for $j\in
\{1,3,\cdots,k-2\}$,

$c_k(v_{i,0}^{(1)})=1$, $c_k(v_{i,0}^{(2)})=k-1$, for $i\in
\{1,2,\cdots,n\}$,

$c_k(u_{i,j}^{(1)})=1$, $c_k(u_{i,j}^{(2)})=k-1$, for
$u_{i,j}^{(1)},u_{i,j}^{(2)} \in V_k^{(0)}$,

$c_k(w_{i,j}^{(1)})=2,$ $c_k(w_{i,j}^{(2)})=3,$ for
$w_{i,j}^{(1)},w_{i,j}^{(2)} \in V_k^{(1)}$,

$c_k(v_{i,k})=c(v_i)$, for $i\in \{1,2,\cdots,n\}$.

{\bf Subcase 2.1.} $v\in V_k^{(p)}$, $w\in V_k^{(q)}$, where $p,q\in
\{1,3,\cdots, k-2\}$.

If $v=v_{i,p}^{(\alpha)}$, $w=v_{j,q}^{(\beta)}$ for
$\alpha,\beta\in \{1,2\},$ then
$vv_{i,p-2}^{(1)}v_{i,p-4}^{(1)}\cdots
v_{i,0}^{(1)}v_{j,0}^{(2)}v_{j,1}^{(2)}\cdots v_{j,q-2}^{(2)}w$ is
the rainbow vertex-connected path between $v$ and $w$.

If $v=v_{i,p}^{(\alpha)}$, $w=w_{i,j}^{(\beta)}$ for
$\alpha,\beta\in \{1,2\},$ then
$vv_{i,p-2}^{(1)}v_{i,p-4}^{(1)}\cdots v_{i,0}^{(1)}u_{i,j}^{(2)}w$
is the rainbow vertex-connected path between $v$ and $w$.

If $v=w_{i_1,j_1}^{(\alpha)}$, $w=w_{i_2,j_2}^{(\beta)}$ for
$\alpha,\beta\in \{1,2\},$ then
$vu_{i_1,j_1}^{(1)}u_{i_2,j_2}^{(2)}w$ is the rainbow
vertex-connected path between $v$ and $w$.

{\bf Subcase 2.2.} $v=v_{i,k}$, $w\in V_k^{(q)}$, where $q\in
\{1,3,\cdots, k-2\}$.

If $w=v_{j,q}^{(\alpha)}$ for $\alpha\in \{1,2\},$ then
$vv_{i,k-2}^{(2)}v_{i,k-4}^{(1)}\cdots v_{i,1}^{(1)}
v_{i,0}^{(1)}v_{j,0}^{(2)}v_{j,2}^{(2)}\cdots v_{j,q-2}^{(2)}w$ is
the rainbow vertex-connected path between $v$ and $w$.

If $w=w_{i,j}^{(\alpha)}$ for $\alpha \in \{1,2\},$ then
$vv_{i,k-2}^{(2)}v_{i,k-4}^{(1)}\cdots
v_{i,1}^{(1)}v_{i,0}^{(1)}u_{i,j}^{(2)}w$ is the rainbow
vertex-connected path between $v$ and $w$.

{\bf Subcase 2.3.} $v=v_{i,k}$, $w=v_{j,k}$.

If $(v_{i,k},v_{j,k})\in P_k$, then $(v_i,v_j)\in P$. By the
vertex-coloring $c$ of $G$, there is a rainbow vertex-connected path
between $v_i$ and $v_j$ in $G$. Correspondingly, since
$c_k(v_{i,k})=c(v_i)$, there is a rainbow vertex-connected path
between $v_{i,k}$ and $v_{j,k}$ in $G_k$.

If $(v_{i,k},v_{j,k})\notin P_k$, then
$v_{i,k}v_{i,k-2}^{(1)}v_{i,k-4}^{(1)}\cdots
v_{i,3}^{(1)}w_{i,j}^{(1)}u_{i,j}^{(1)}w_{i,j}^{(2)}v_{j,3}^{(2)}\cdots
v_{j,k-2}^{(2)}v_{j,k}$ is the rainbow vertex-connected path between
$v_{i,k}$ and $v_{j,k}$. \qed

\noindent {\bf Proof of Theorem \ref{theorem2}:} From the above
Lemmas 1 and 2, the first part of Theorem \ref{theorem2}, the
NP-Hardness, follows immediately.

In the following we will prove the second part of Theorem
\ref{theorem2}. Recall that a problem belongs to NP-class if given
any instance of the problem whose answer is ``yes", there is a
certificate validating this fact which can be checked in polynomial
time. For any fixed integer $k$, to prove the problem of deciding
whether $rvc(G)\leq k$ is in NP-class, we can choose a rainbow
$k$-vertex-coloring of $G$ as a certificate. For checking a rainbow
$k$-vertex-coloring, we only need to check that $k$ colors are used
and for any two vertices $u$ and $v$ of $G$, there exists a rainbow
vertex-connected path between $u$ and $v$. Notice that for any two
vertices $u$ and $v$ of $G$, there are at most $n^{\ell-1}$ $u-v$
paths of length $\ell$, since if we let $P = uv_1v_2\cdots
v_{\ell-1}v$, then there are less than $n$ choices for each $v_i$
$(i \in \{1, 2,\ldots, \ell-1\})$. Therefore, $G$ contains at most
$\sum_{\ell=1}^{k+1}n^{\ell-1}= \frac{n^{k+1}-1}{n}\leq n^k$ $u-v$
paths of length at most $k+1$. Then, check these paths in turn until
one finds one path whose internal vertices have distinct colors. It
follows that the time used for checking is at most $O(n^k \cdot n^2
\cdot n^2) = O(n^{k+4})$. Since $k$ is a fixed integer, we conclude
that the certificate can be checked in polynomial time, which
implies that the problem of deciding whether $rvc(G)\leq k$ belongs
to NP-class, and therefore it is NP-Complete. \qed


\begin{thebibliography}{1}

\bibitem{AN}
P. Ananth, M. Nasre, New hardness results in rainbow connectivity,
arXiv:1104.2074v1 [cs.CC] 2011.

\bibitem{BM} J.A. Bondy, U.S.R. Murty, Graph Theory, GTM 244, Springer, 2008.

\bibitem{CLRTY}
Y. Caro, A. Lev, Y. Roditty, Z. Tuza, R. Yuster,  On rainbow
connection, {\it Electron J. Combin.} {\bf 15}(2008), R57.

\bibitem{CFMY}
S. Chakraborty, E. Fischer, A. Matsliah, R. Yuster, Hardness and
algorithms for rainbow connectivity, {\it J. Comb. Optim.} {\bf
21}(2011), 330--347.

\bibitem{CJMZ}
G. Chartrand, G.L. Johns, K.A. McKeon, P. Zhang, Rainbow connection
in graphs, {\it Math. Bohemica} {\bf 133}(2008), 85--98.

\bibitem{CLS}
L. Chen, X. Li, Y. Shi, The complexity of determining the rainbow
vertex-connection of graphs, {\it Theoretical Computer Science} {\bf
412}(2011), 4531--4535.

\bibitem{KY}
M. Krivelevich, R. Yuster, The rainbow connection of a graph is (at
most) reciprocal to its minimum degree, {\it J. Graph Theory} {\bf
63}(2010), 185--191.

\bibitem{SX}
S. Li, X. Li, Note on the complexity of deciding the rainbow
connectedness for bipartite graphs, arXiv:1109.5534v1 [math.CO]
2011.

\bibitem{LS}
X. Li, Y. Shi, On the rainbow vertex-connection, arXiv:1012.3504v1
[math.CO] 2010.

\bibitem{S}
I. Schiermeyer, Rainbow connection in graphs with minimum degree
three, IWOCA 2009, {\it LNCS} {\bf 5874}(2009), 432--437.

\end{thebibliography}
\end{document}